\documentclass[a4paper,11pt]{article}
\usepackage{CARen}
\PassOptionsToPackage{hyphens}{url}
\usepackage{hyperref}

\newcommand{\RR}{\mathbb R}

\begin{document}

\setcounter{footnote}{0}
\setcounter{figure}{0}


\Abschnitt
{}
{}
{}

\vspace{3mm}

\Aufsatz
{Hands-on Tropical Geometry}
{Hands-on Tropical Geometry}
{H. Gangl, Y. Ren, Z. Urbancic}
{NameAuthor}
{H. Gangl, Durham University\\ Y. Ren, Durham University\\ Z. Urbancic, Durham University}
{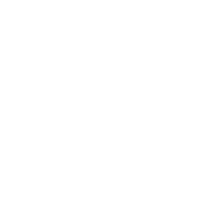}
{{herbert.gangl@durham.ac.uk, yue.ren2@durham.ac.uk, ziva.urbancic@durham.ac.uk}}


\vspace{3mm}
\begin{multicols}{2}
\noindent


\Ueberschrift{Introduction\vphantom{y}}{sec:introduction}

Tropical varieties are piecewise linear structures that arise in many areas inside and outside mathematics. They describe the loci of indifference prices in product mix auctions \cite{BaldwinKlemperer2019}, spaces of phylogenetic trees \cite{SpeyerSturmfels04,MLYK2018}, or images of solutions of polynomial systems under component-wise valuation \cite{MaclaganSturmfels2015}. Geometrically, tropical varieties are strongly tied to their classical algebro-geometric counterparts, they are strikingly similar in some aspects yet also fascinatingly different in others. However, while there has been a significant effort in modelling algebraic surfaces and curves - from Schilling's white plaster casts popularized by Klein in the 19th century~\cite{Schilling1911,HalverscheidLabs2019} to 3D-printed models of today \cite{MathSculptures} \cite{BertiniReal} - nothing comparable exists for tropical varieties despite the simplicity of their polyhedral nature.

This article aims to address this gap by presenting a comprehensive guide on how to create 3D-printable models of tropical surfaces, tropical curves, and combinations thereof. We will make use the following open-source software, which are freely available on their respective websites for all platforms:
\begin{description}[leftmargin=*]
\item[\textsc{polymake} \cite{Polymake}] (\url{https://polymake.org})\\ for creating mathematically accurate tropical models. \vspace{-1mm}
\item[\textsc{OpenSCAD} \cite{OpenSCAD}] (\url{https://openscad.org})\\ for solidifying the tropical models. 
\end{description}
We will only focus on the creation of 3D models in \textsc{OpenSCAD}, which in turn is capable of exporting to a variety of file formats common in 3D printing such as \texttt{STL}. Exporting color is also supported by 3rd party scripts such as \textsc{ColorSCAD}. The actual 3D printing process depends significantly on the printer, the material, and the associated printing soft\-ware used. 

All scripts and templates in this report can be found in the official polymake repository and will be part of the official polymake distribution from version 4.6 onward. The notebook is also publicly available under
\begin{center}\footnotesize
  \texttt{https://polymake.org/doku.php/user\_guide}\\\texttt{/tutorials/master/hands\_on\_tropical\_geometry}
\end{center}

\Ueberschrift{Summary}{sec:summary}
This article assumes some familiarity with the basic objects in tropical geometry. If any terminology is unkown, refer to \cite{MaclaganSturmfels2015} and \cite{Joswig2021}.

The starting point for creating a 3D-model is one or more polyhedral complexes in $\RR^3$. In \textsc{polymake}, objects of type \texttt{fan::PolyhedralComplex} can always be constructed manually by specifying their \texttt{POINTS} and \texttt{INPUT\_POLYTOPES}. In special cases, which are covered in the upcoming sections, they can also be constructed via some shortcuts.

Given polyhedral complexes in \textsc{polymake}, the 3D-model can then be created in four steps:
\begin{enumerate}
\item Construct a bounding box in \textsc{polymake}. This can be done automatically using our script (which draws a box with a prescribed margin around all vertices) or manually by specifying any bounded polytope (which need not be a cuboid). Use ``\texttt{->VISUAL}'' to verify that the bounded box cuts off the polyhedral complex as desired.\vspace{-1.5mm}
\item Export the model from \textsc{polymake} to \textsc{OpenSCAD} using our script.\vspace{-1.5mm}
\item Adjust the thicknesses in \textsc{OpenSCAD}. Use the preview window to verify that the thickness is as desired.\vspace{-1.5mm}
\item Export the model from \textsc{OpenSCAD} to any 3D-printable file format.
\end{enumerate}
In the following sections, we will specifically discuss the cases of
\begin{itemize}
\item a single tropical surface,\vspace{-1.5mm}
\item a single tropical curve,\vspace{-1.5mm}
\item a tropical surface containing a tropical curve.
\end{itemize}
The methods used can however be combined to cover arbitrary arrangements of surfaces and curves.

\pagebreak

\Ueberschrift{Modelling Tropical Surfaces}{sec:tropicalSurfaces}
For an example of producing a 3D-printable model of a tropical surface see the file
\begin{center}
  \texttt{3d\_printing\_template\_surface.pl}.
\end{center}
To run the example, simply copy its contents into any \textsc{polymake} session.
In this section, we briefly explain some of its contents. Note that some variables are renamed in this article due to spacing.

\Ueberschriftu{Constructing a tropical surface in \textsc{polymake}}
An easy way to construct a tropical surface in $\RR^3$ of type \texttt{fan::PolyhedralComplex} in \textsc{polymake} is by constructing and converting a \texttt{tropical::Hypersurface}.
This can be done either by specifying a tropical polynomial as a string, or by specifying an exponent matrix and a coefficient vector:

\begin{lstlisting}[basicstyle=\scriptsize\ttfamily,
  showlines=true,
  showstringspaces=false,
  language=perl,
  keywordstyle=\color{blue}\ttfamily,
  stringstyle=\color{red}\ttfamily,
  commentstyle=\color{green!30!black},
  morekeywords={Hypersurface,Min}]
# Tropical quadratic surface via polynomial:
$f = toTropicalPolynomial("min(1+2*w,1+2*x,1+2*y,
                               1+2*z,w+x,w+y,w+z,
                               x+y,x+z,y+z)");
$T = new tropical::Hypersurface<Min>(POLYNOMIAL=>
                                       $f);
# Tropical quadratic surface via exponent matrix
#   and coefficient vector:
$m = [[2,0,0,0], [1,1,0,0], [1,0,1,0], [1,0,0,1],
      [0,1,1,0], [0,1,0,1], [0,0,1,1], [0,2,0,0],
      [0,0,2,0], [0,0,0,2]];
$c = [1,0,0,0,0,0,0,1,1,1];
$T = new Hypersurface<Min>(MONOMIALS=>$m,
                           COEFFICIENTS=>$c);
# Converting Hypersurface to PolyhedralComplex
$T = new fan::PolyhedralComplex(
       VERTICES=>$T->VERTICES->minor(All,~[1]),
       MAXIMAL_POLYTOPES=>$T->MAXIMAL_POLYTOPES);
# Visualization:
$T->VISUAL;
\end{lstlisting}

\textbf{Important:} Since \textsc{polymake} uses homogeneous coordinates, i.e., it identifies polyhedra in $\RR^3$ with cones in $\RR^4$, tropical polynomials need to be tetravariate and homogeneous instead of trivariate and inhomogeneous.


\begin{figurehere}
 \centering
 \includegraphics[width=0.5\columnwidth]{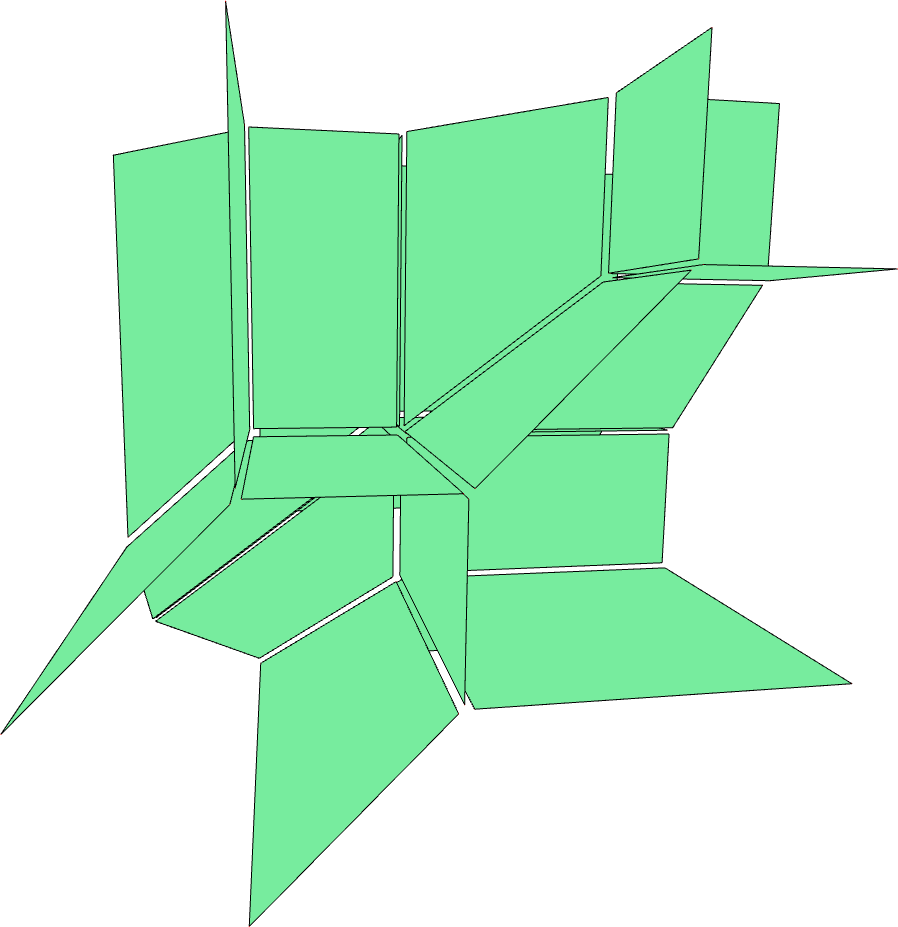}
 \caption{Tropical quadratic surface in \textsc{polymake}.\label{fig:polymakeSurface}}
\end{figurehere}

\Ueberschriftu{Bounding a tropical surface in \textsc{polymake} and exporting it to \textsc{OpenSCAD}}
Constructing a bounding box can be done using the command \texttt{generateBoundingBox} or by specifying a custom \texttt{polytope}, which need not be a cuboid:

\begin{lstlisting}[basicstyle=\scriptsize\ttfamily,
  showlines=true,
  showstringspaces=false,
  language=perl,
  keywordstyle=\color{blue}\ttfamily,
  stringstyle=\color{red}\ttfamily,
  commentstyle=\color{green!30!black},
  morekeywords={Hypersurface,Min}]
# Constructing bounding box automatically:
$bBox = generateBoundingBox($T);
# Constructing bounding box manually:
$bBox = scale(cube(3),2);
\end{lstlisting}

Note: \texttt{generateBoundingBox(\$T)} draws a box around all vertices of \texttt{\$T} with a margin of $1$, \texttt{generateBoundingBox(\$T,3,4,5)} draws a box around all vertices of \texttt{\$T} with margins $3$, $4$, $5$ in $x$-, $y$-, and $z$-direction, respectively.

The command \texttt{intersectWithBoundingBox} intersects the \texttt{fan::PolyhedralComplex} with the bounding box, \texttt{->VISUAL} visualizes their intersection:

\begin{lstlisting}[basicstyle=\scriptsize\ttfamily,
  showlines=true,
  showstringspaces=false,
  language=perl,
  keywordstyle=\color{blue}\ttfamily,
  stringstyle=\color{red}\ttfamily,
  commentstyle=\color{green!30!black},
  morekeywords={Hypersurface,Min}]
# Constructing bounded PolyhedralComplex:
$Tbounded = intersectWithBoundingBox($T,$bBox);
$Tbounded->VISUAL;
\end{lstlisting}

\medskip
The command \texttt{generateSCADFileForSurface} exports the tropical surface to \textsc{OpenSCAD}. It requires a bounded polyhedral complex and a filename. If the file already exists, it will be overwritten:

\begin{lstlisting}[basicstyle=\scriptsize\ttfamily,
  showlines=true,
  showstringspaces=false,
  language=perl,
  keywordstyle=\color{blue}\ttfamily,
  stringstyle=\color{red}\ttfamily,
  commentstyle=\color{green!30!black},
  morekeywords={Hypersurface,Min}]
$filename = "foo.scad";
generateSCADFileForSurface($Tbounded,$filename);
\end{lstlisting}

\Ueberschriftu{Solidifying a tropical surface in \textsc{OpenSCAD} and exporting it for 3D-printing}
To solidify the surface into a three-dimensional model, open the exported file in \textsc{OpenSCAD}. See Figure~\ref{fig:openscadSurface} for a preview of the model.

\begin{figurehere}
 \centering
 \includegraphics[height=37mm]{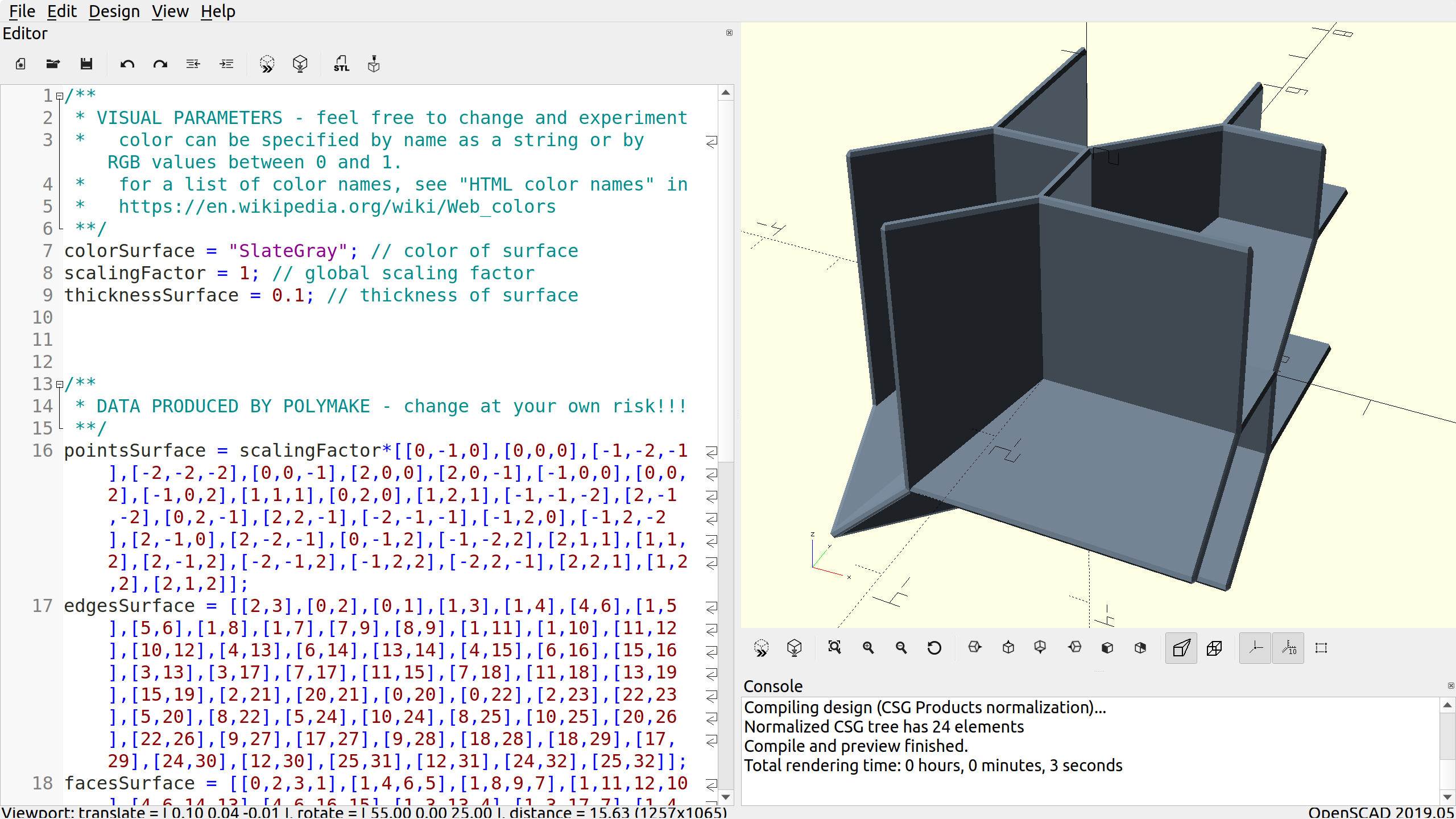}
 \caption{Solidified tropical quadratic surface in \textsc{OpenSCAD}.\label{fig:openscadSurface}}
\end{figurehere}

The top of the file contains parameters which control the thickness of the model amongst other things:

\begin{lstlisting}[basicstyle=\scriptsize\ttfamily,
  showlines=true,
  showstringspaces=false,
  language=perl,
  keywordstyle=\color{blue}\ttfamily,
  stringstyle=\color{red}\ttfamily,
  commentstyle=\color{green!30!black},
  morekeywords={Hypersurface,Min}]
colorSurface = "SlateGray"; // color of surface
scalingFactor = 1; // global scaling factor
thicknessSurface = 0.05; // thickness of surface
\end{lstlisting}

\begin{itemize}
\item \texttt{colorSurface} is the color of the surface in the render. It can be specified either by an HTML color name or by RGB value.
\item \texttt{scalingFactor} is a factor by which the polyhedral complex is scaled. It can be used to scale the model to a desired size.
\item \texttt{thicknessSurface} controls the thickness of the surface. \textsc{OpenSCAD} solidifies the surface by replacing every vertex with a ball with specified thickness and taking convex hulls.
\end{itemize}

\textbf{Important:} In order to 3D-print the tropical surface, every face needs to have a minimal thickness. As a rule of thumb for printing with PLA, we recommend a minimal thickness of $2$mm for a print of size $(10\text{cm})^3$.

\Ueberschrift{Tropical Curves}{sec:tropicalCurves}

For an example of producing a 3D-printable model of a tropical curve see the file
\begin{center}
  \texttt{3d\_printing\_template\_curve.pl}.
\end{center}
To run the example, simply copy its contents into any \textsc{polymake} session.
In this section, we briefly explain some of its contents. Note that some variables are renamed due to spacing.

\Ueberschriftu{Constructing a tropical curve in \textsc{polymake}}
An easy way to construct a tropical curve in $\RR^3$ of type \texttt{fan::PolyhedralComplex} in \textsc{polymake} is by intersecting two tropical surfaces. Note however that not all tropical curves can be constructed this way.

\begin{lstlisting}[basicstyle=\scriptsize\ttfamily,
  showlines=true,
  showstringspaces=false,
  language=perl,
  keywordstyle=\color{blue}\ttfamily,
  stringstyle=\color{red}\ttfamily,
  commentstyle=\color{green!30!black},
  morekeywords={Hypersurface,Min}]
# Tropical sextic curve that is the intersection
# of a tropical quadratic and cubic surface
$Mq = [[2,0,0,0], [0,2,0,0], [0,0,2,0], [0,0,0,2],
       [1,1,0,0], [1,0,1,0], [1,0,0,1],
       [0,1,1,0], [0,1,0,1], [0,0,1,1]];
$Cq = [1, -1/4, -2/4, -3/4, -3/4,
       -4/4, -5/4, 2/4, 0, -2/4];
$Tq = new tropical::Hypersurface<Min>(
                    MONOMIALS=>$Mq,
                    COEFFICIENTS=>$Cq);
$Mc = [[3,0,0,0], [0,3,0,0], [0,0,3,0], [0,0,0,3],
       [1,1,1,0], [1,1,0,1], [1,0,1,1], [0,1,1,1],
       [2,1,0,0], [2,0,1,0], [2,0,0,1],
       [1,2,0,0], [1,0,2,0], [1,0,0,2],
       [0,2,1,0], [0,2,0,1], [0,1,2,0],
       [0,1,0,2], [0,0,2,1], [0,0,1,2]];
$Cc = [3,3,3,3,0,0,0,0,1,1,1,1,1,1,1,1,1,1,1,1];
$Tc = new tropical::Hypersurface<Min>(
                MONOMIALS=>$Mc,
                COEFFICIENTS=>$Cc);
$Ts = tropical::intersect($Tq,$Tc);
# Converting Cycle to PolyhedralComplex
$Ts = new fan::PolyhedralComplex(
       VERTICES=>$Ts->VERTICES->minor(All,~[1]),
       MAXIMAL_POLYTOPES=>$Ts->MAXIMAL_POLYTOPES);
\end{lstlisting}
\textbf{Important:} \texttt{intersect} is the intersection of tropical cycles. In terms of polyhedral complexes, it therefore computes the stable intersection, not the set-theoretic intersection.

Alternatively and in particular for tropical curves which are not realizable, they can also be manually constructed as a one-dimensional polyhedral complex by specifying vertices and maximal polyhedra:

\begin{lstlisting}[basicstyle=\scriptsize\ttfamily,
  showlines=true,
  showstringspaces=false,
  language=perl,
  keywordstyle=\color{blue}\ttfamily,
  stringstyle=\color{red}\ttfamily,
  commentstyle=\color{green!30!black},
  morekeywords={Hypersurface,Min}]
# constructing a tropical genus-2 cubic curve
$V = [[1,-1,-1,0],  [1,-4,-3,0], [1,-8,-5,0],
      [1,-9,-6,0],  [1,-9,-7,0], [1,-8,-7,0],
      [1,-10,-8,0], [1,1,1,0],   [1,4,4,1],
      [1,8,8,1],    [1,9,9,2],   [1,9,9,3],
      [1,8,8,3],    [1,10,10,2], [0,-1,0,0],
      [0,0,-1,0],   [0,0,0,-1],  [0,1,1,1]];
$E = [[0,1],   [1,2],  [2,3],   [3,4],   [4,5],
      [5,1],   [4,6],  [0,15],  [5,15],  [6,15],
      [6,14],  [3,14], [2,14],  [0,7],   [7,8],
      [8,9],   [9,10], [10,11], [11,12], [12,8],
      [10,13], [7,16], [9,16],  [13,16], [13,17],
      [11,17], [12,17]];
$Tc = new fan::PolyhedralComplex(VERTICES=>$V,MAXIMAL_POLYTOPES=>$E);
\end{lstlisting}


\Ueberschriftu{Framing a tropical curve in \textsc{polymake} and exporting it to \textsc{OpenSCAD}}
Constructing a bounding box can be done using the command \texttt{generateBoundingBox} as before:

\begin{lstlisting}[basicstyle=\scriptsize\ttfamily,
  showlines=true,
  showstringspaces=false,
  language=perl,
  keywordstyle=\color{blue}\ttfamily,
  stringstyle=\color{red}\ttfamily,
  morekeywords={Hypersurface,Min}]
$bBox = generateBoundingBox($Ts);
$Tbounded = intersectWithBoundingBox($Ts,$bBox);
\end{lstlisting}

For the sake of stability, we recommend constructing a frame for the tropical curve. This is best done using a tropical surface containing it:
\begin{lstlisting}[basicstyle=\scriptsize\ttfamily,
  showlines=true,
  showstringspaces=false,
  language=perl,
  keywordstyle=\color{blue}\ttfamily,
  stringstyle=\color{red}\ttfamily,
  morekeywords={Hypersurface,Min}]
$Tq = new fan::PolyhedralComplex(
       VERTICES=>$Tq->VERTICES->minor(All,~[1]),
       MAXIMAL_POLYTOPES=>$Tq->MAXIMAL_POLYTOPES);
$Tframe = intersectWithBoundingBoxForFraming($Tq,$bBox);
\end{lstlisting}

Exporting the frame and the tropical curve to \textsc{OpenSCAD} can be done using the command \texttt{generateSCADFileForCurve}. It requires two bounded polyhedral complexes and a filename. If the file already exists, it will be overwritten:

\begin{lstlisting}[basicstyle=\scriptsize\ttfamily,
  showlines=true,
  showstringspaces=false,
  language=perl,
  keywordstyle=\color{blue}\ttfamily,
  stringstyle=\color{red}\ttfamily,
  commentstyle=\color{green!30!black},
  morekeywords={Hypersurface,Min}]
$filename = "foo.scad";
generateSCADFileForCurve($Tframe,$Ts,$filename);
\end{lstlisting}

\Ueberschriftu{Solidifying a tropical curve in \textsc{OpenSCAD} and exporting it for 3D-printing}
To solidify the curve into a three-dimensional model, open the exported file in \textsc{OpenSCAD}. See Figures~\ref{fig:openscadCurve} and \ref{fig:openscadCurve2} for a preview of the model.

\begin{figurehere}
 \centering
 \includegraphics[height=37mm]{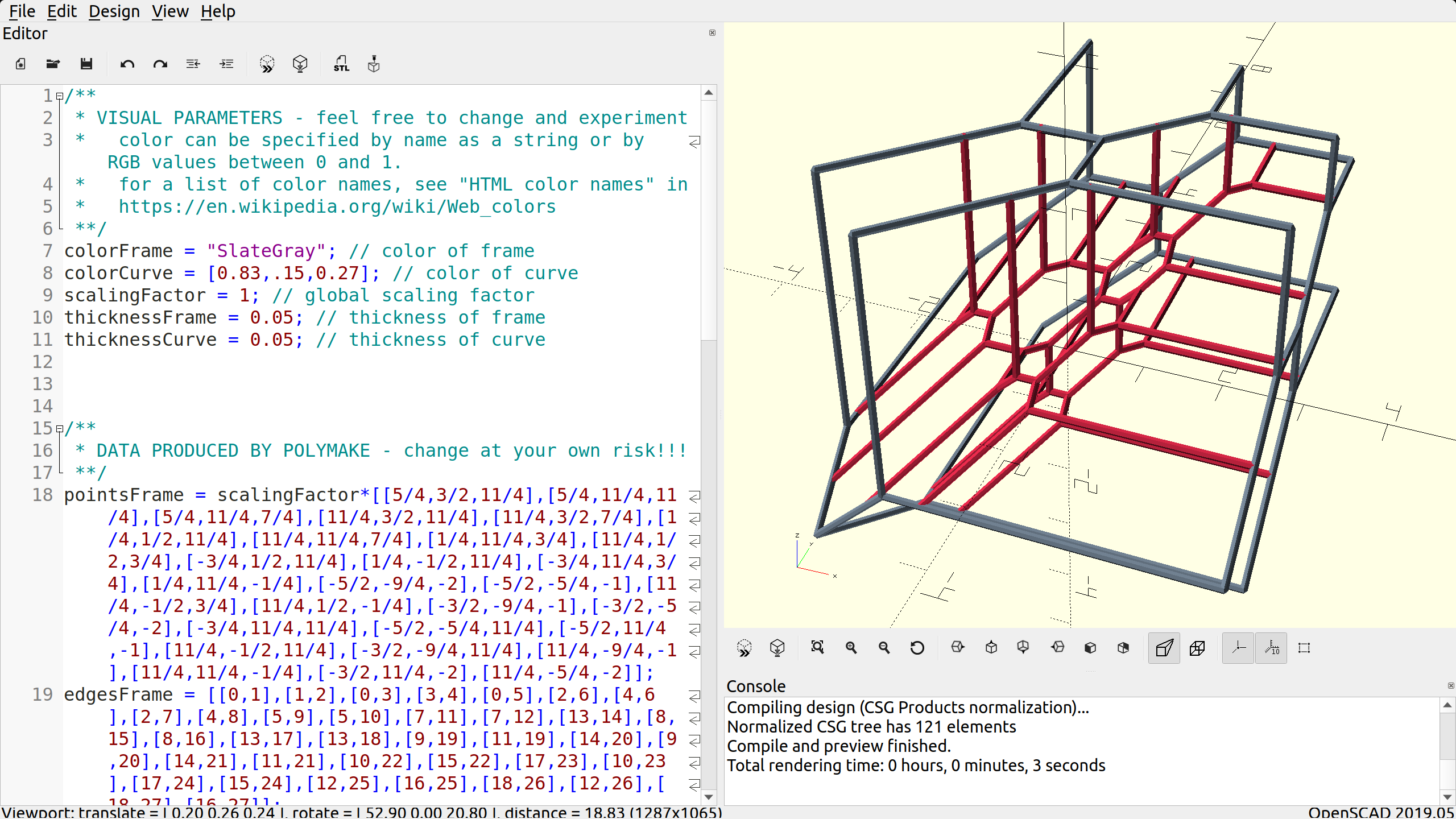}
 \caption{Framed and solidified tropical sextic curve \textsc{OpenSCAD}.\label{fig:openscadCurve}}
\end{figurehere}

\begin{figurehere}
 \centering
 \includegraphics[height=37mm]{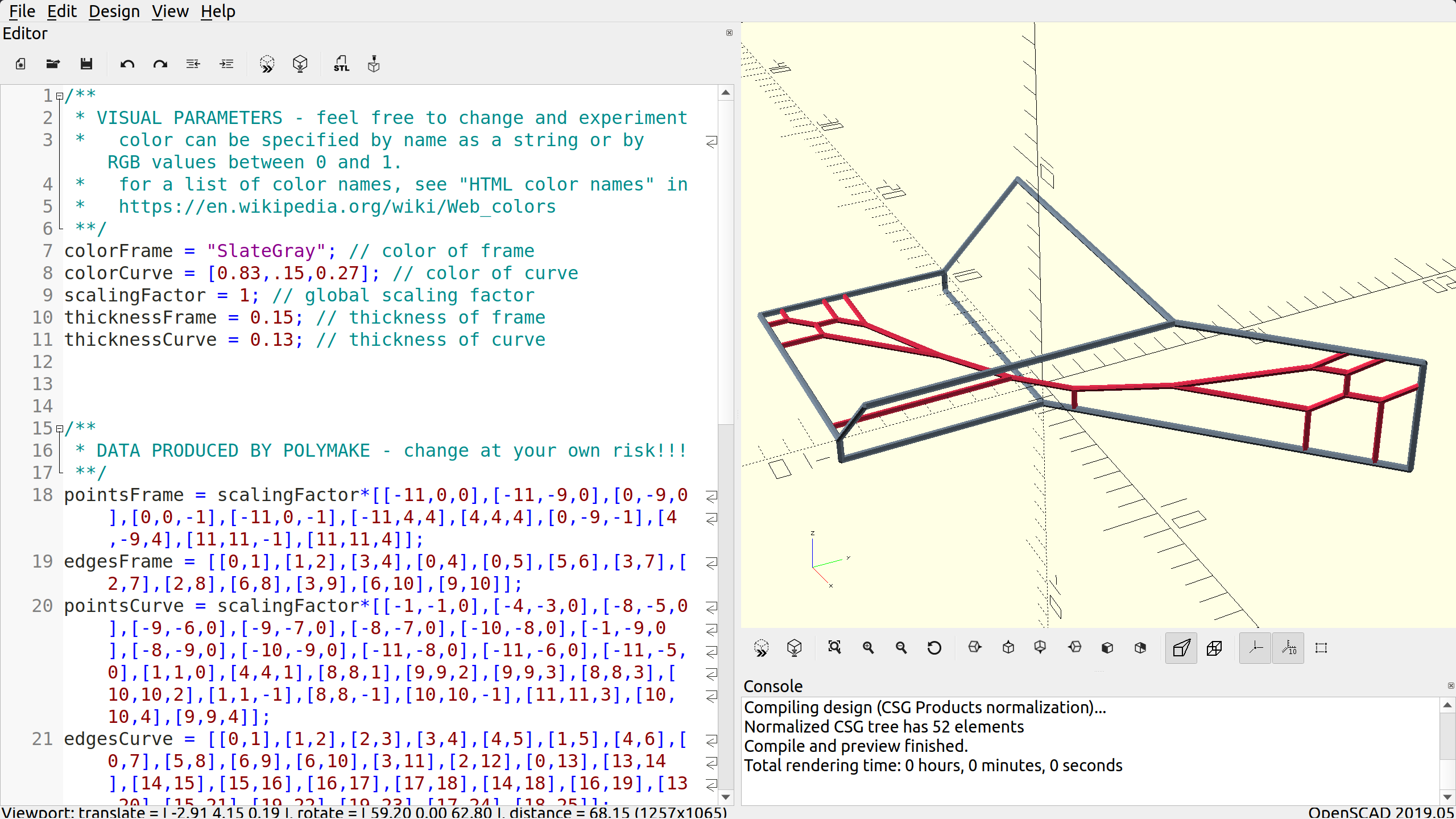}
 \caption{Framed and solidified genus-2 tropical cubic curve in \textsc{OpenSCAD}.\label{fig:openscadCurve2}}
\end{figurehere}

As before, the top of the file contains parameters which control the thickness of the model amongst other things:

\begin{lstlisting}[basicstyle=\scriptsize\ttfamily,
  showlines=true,
  showstringspaces=false,
  language=perl,
  keywordstyle=\color{blue}\ttfamily,
  stringstyle=\color{red}\ttfamily,
  commentstyle=\color{green!30!black},
  morekeywords={Hypersurface,Min}]
colorFrame = "SlateGray"; // color of frame
colorCurve = [0.83,.15,0.27]; // color of curve
scalingFactor = 1; // global scaling factor
thicknessFrame = 0.05; // thickness of frame
thicknessCurve = 0.05; // thickness of curve
\end{lstlisting}

\begin{itemize}
\item \texttt{colorFrame} and \texttt{colorCurve} are the colors of the surface in the render. They can be specified either by an HTML color name or by RGB value.
\item \texttt{scalingFactor} is a factor by which the polyhedral complex is scaled. It can be used to scale the model to the desired size.
\item \texttt{thicknessFrame} and \texttt{thicknessCurve} control the thickness of the frame and curve. \textsc{OpenSCAD} solidifies both by replacing every vertex with a ball with diameter equal to the specified radius and taking convex hulls.
\end{itemize}

\Ueberschrift{Modelling Tropical Curves on\\ Tropical Surfaces}{sec:tropicalCurvesSurfaces}

For an example of producing a 3D-printable model of a tropical curve on a tropical surface see the file
\begin{center}\footnotesize
  \texttt{3d\_printing\_template\_surface\_and\_curve.pl}.
\end{center}
To run the example, simply copy its contents into any \textsc{polymake} session.
In this section, we briefly explain some of its contents. Note that some variables are renamed due to spacing.

\Ueberschriftu{Constructing a curve and surface in \textsc{polymake}}
The previous section showed an easy way to construct a tropical curve lying on a tropical surface in \textsc{polymake}: start with the surface and construct the curve by intersecting it with another surface. However, as mentioned before, this may not always possible or easy to do. One option that always works is to construct the curve manually by specifying its vertices and edges:

\begin{lstlisting}[basicstyle=\scriptsize\ttfamily,
  showlines=true,
  showstringspaces=false,
  language=perl,
  keywordstyle=\color{blue}\ttfamily,
  stringstyle=\color{red}\ttfamily,
  commentstyle=\color{green!30!black},
  morekeywords={Hypersurface,Min}]
# constructing tropical plane via polynomial
$l=toTropicalPolynomial("max(y,z,w)",qw(w x y z));
$Tl = new tropical::Hypersurface<Max>(
        POLYNOMIAL=>$l);
$Tl = new fan::PolyhedralComplex(
        POINTS=>$Tl->VERTICES->minor(All,~[1]),
        INPUT_POLYTOPES=>$Tl->MAXIMAL_POLYTOPES,
        INPUT_LINEALITY=>$Tl->LINEALITY_SPACE->minor(All,~[1]));

# constructing tropical curve manually:
$V = [[1,0,0,0],[1,-8,0,0],[1,-8,-6,0],
      [1,12,0,0],[0,0,-1,0],[1,-6,-6,0],
      [1,-11,-3,0],[1,-11,-6,0],[0,-1,0,0],
      [1,-12,-7,0],[1,4,0,-2],[1,10,0,-2],
      [1,8,0,-6],[0,0,0,-1],[0,1,1,1],
      [1,-4,6,6],[1,-4,5,5],[1,-3,5,5],
      [1,-5,4,4],[1,-3,3,3],[1,-5,3,3]];
$E = [[3,4],[0,5],[2,5],[4,5],[1,6],[6,7],
      [6,8],[2,7],[7,9],[8,9],[4,9],[10,11],
      [3,11],[11,12],[0,10],[10,12],[12,13],
      [1,13],[14,15],[8,15],[15,16],[14,17],
      [16,17],[16,18],[17,19],[8,18],[18,20],
      [0,19],[19,20],[1,20],[3,14]];
$Tq = new fan::PolyhedralComplex(POINTS=>$V,
        INPUT_POLYTOPES=>$E);
\end{lstlisting}

\begin{figurehere}
 \centering
 \includegraphics[width=0.8\columnwidth]{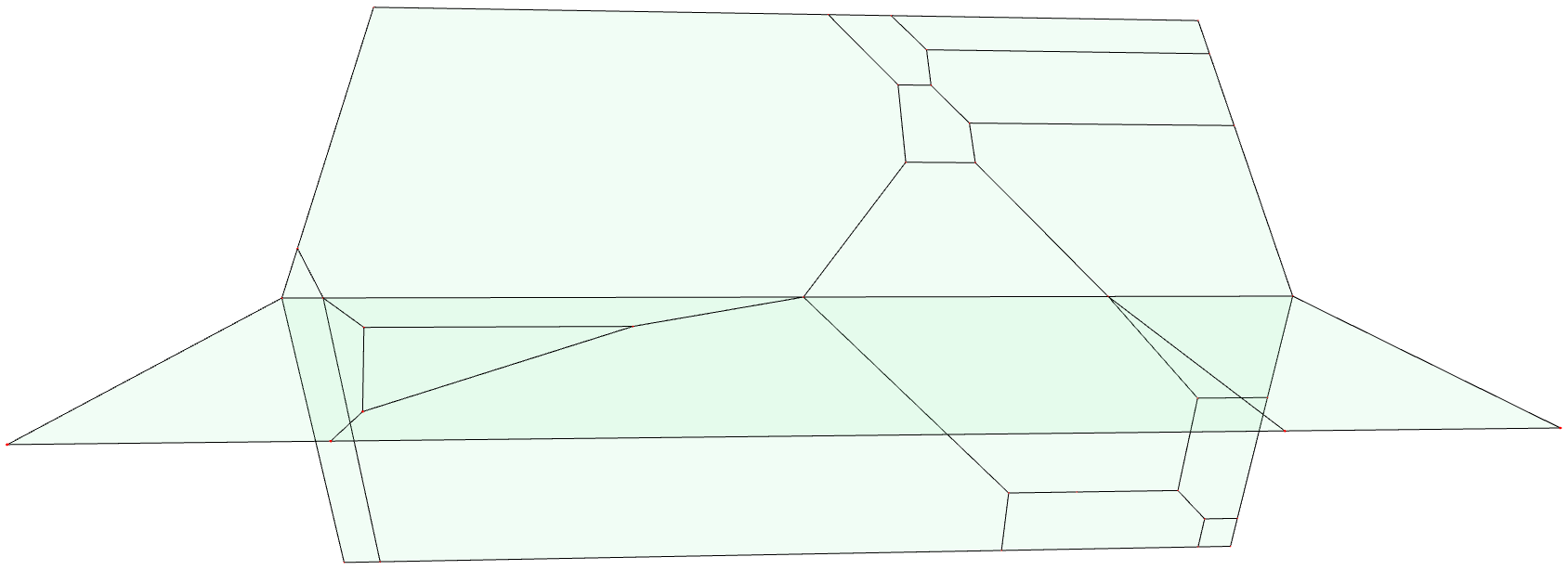}
 \caption{Tropical plane with a tropical quartic curve in \textsc{polymake}.\label{fig:polymakeSurface}}
\end{figurehere}

\Ueberschriftu{Bounding a tropical surface in \textsc{polymake} and exporting it to \textsc{OpenSCAD}}
Constructing a bounding box is easiest done using the command \texttt{generateBoundingBox} as before:

\begin{lstlisting}[basicstyle=\scriptsize\ttfamily,
  showlines=true,
  showstringspaces=false,
  language=perl,
  keywordstyle=\color{blue}\ttfamily,
  stringstyle=\color{red}\ttfamily,
  commentstyle=\color{green!30!black},
  morekeywords={Hypersurface,Min}]
$bBox = generateBoundingBox($Tq);
$TlBounded = intersectWithBoundingBox($Tl,$bBox);
$TqBounded = intersectWithBoundingBox($Tq,$bBox);
\end{lstlisting}

\medskip
Exporting the tropical surface and curve to \textsc{OpenSCAD} can be done using the command \texttt{generateSCADFileForSurfaceAndCurve}. It requires two bounded polyhedral complexes and a filename. If the file already exists, it will be overwritten:

\begin{lstlisting}[basicstyle=\scriptsize\ttfamily,
  showlines=true,
  showstringspaces=false,
  language=perl,
  keywordstyle=\color{blue}\ttfamily,
  stringstyle=\color{red}\ttfamily,
  commentstyle=\color{green!30!black},
  morekeywords={Hypersurface,Min}]
$filename = "foo.scad";
generateSCADFileForSurfaceAndCurve($TlBounded,$TqBounded,$filename);
\end{lstlisting}

\Ueberschriftu{Solidifying a tropical surface with a tropical curve in \textsc{OpenSCAD} and exporting them for 3D-printing}
To solidify the tropical surface with a tropical curve into a three-dimensional model, open the exported file in \textsc{OpenSCAD}. See Figure~\ref{fig:openscadSurfaceAndCurve} for a preview of the model.

\begin{figurehere}
 \centering
 \includegraphics[height=37mm]{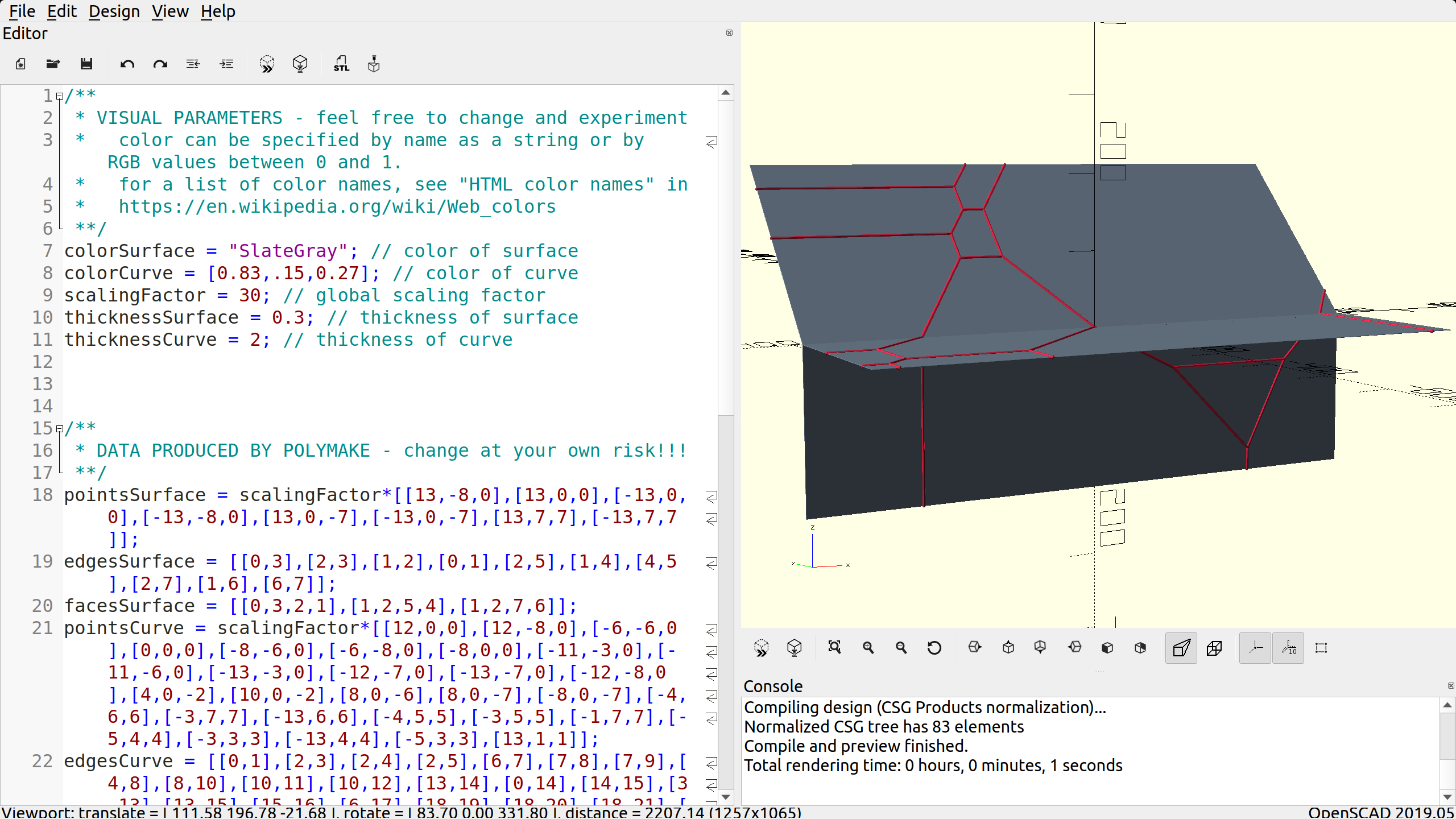}
 \caption{Solidified tropical plane with a tropical quartic curve in \textsc{OpenSCAD}.\label{fig:openscadSurfaceAndCurve}}
\end{figurehere}

The top of the file contains parameters which control the thickness of the model amongst other things:

\begin{lstlisting}[basicstyle=\scriptsize\ttfamily,
  showlines=true,
  showstringspaces=false,
  language=perl,
  keywordstyle=\color{blue}\ttfamily,
  stringstyle=\color{red}\ttfamily,
  commentstyle=\color{green!30!black},
  morekeywords={Hypersurface,Min}]
colorSurface = "SlateGray"; // color of surface
colorCurve = [0.83,.15,0.27]; // color of curve
scalingFactor = 1; // global scaling factor
thicknessSurface = 0.01; // thickness of surface
thicknessCurve = 0.1; // thickness of curve
\end{lstlisting}

\begin{itemize}
\item \texttt{colorSurface} and \texttt{colorCurve} are the colors of the surface and curve in the render. It can be specified either by an HTML color name or by RGB value.
\item \texttt{scalingFactor} is a factor by which the polyhedral complex is scaled. It can be used to scale the model to the desired size.
\item \texttt{thicknessSurface} and \texttt{thicknessCurve} control the thickness of the surface and of the curve, respectively. \textsc{OpenSCAD} solidifies the surface by replacing every vertex with a ball with diameter equal to the specified radius and taking convex hulls.
\end{itemize}



\Ueberschrift{Gallery}{sec:gallery}

\begin{figurehere}
 \centering
 \includegraphics[width=0.8\columnwidth]{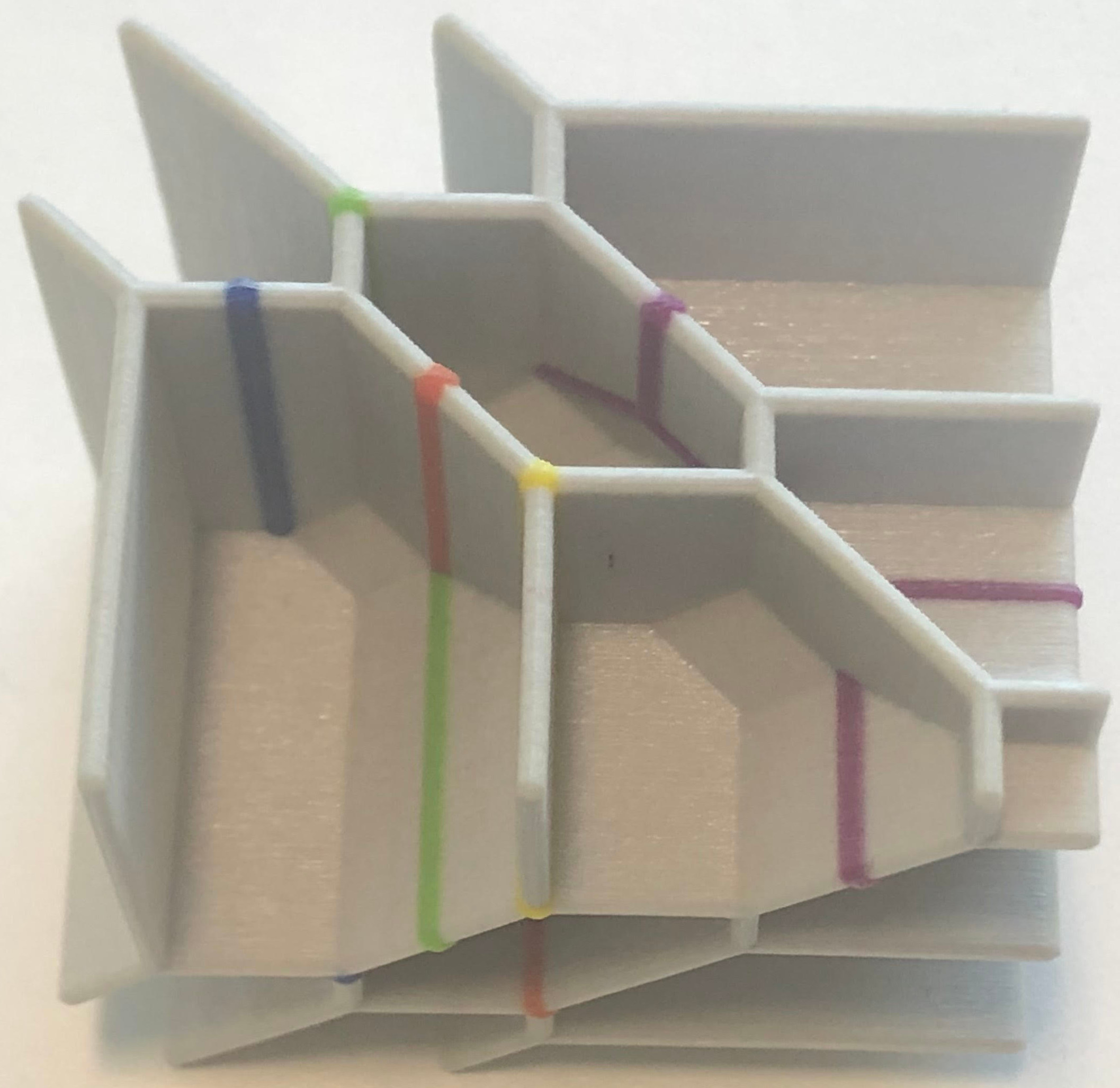}
 \caption{A tropical cubic surface from \cite{PV22} with 5 of 27 tropical lines drawn (one per motif).\label{fig:cubicSurface}}
\end{figurehere}

\begin{figurehere}
 \centering
 \includegraphics[height=37mm]{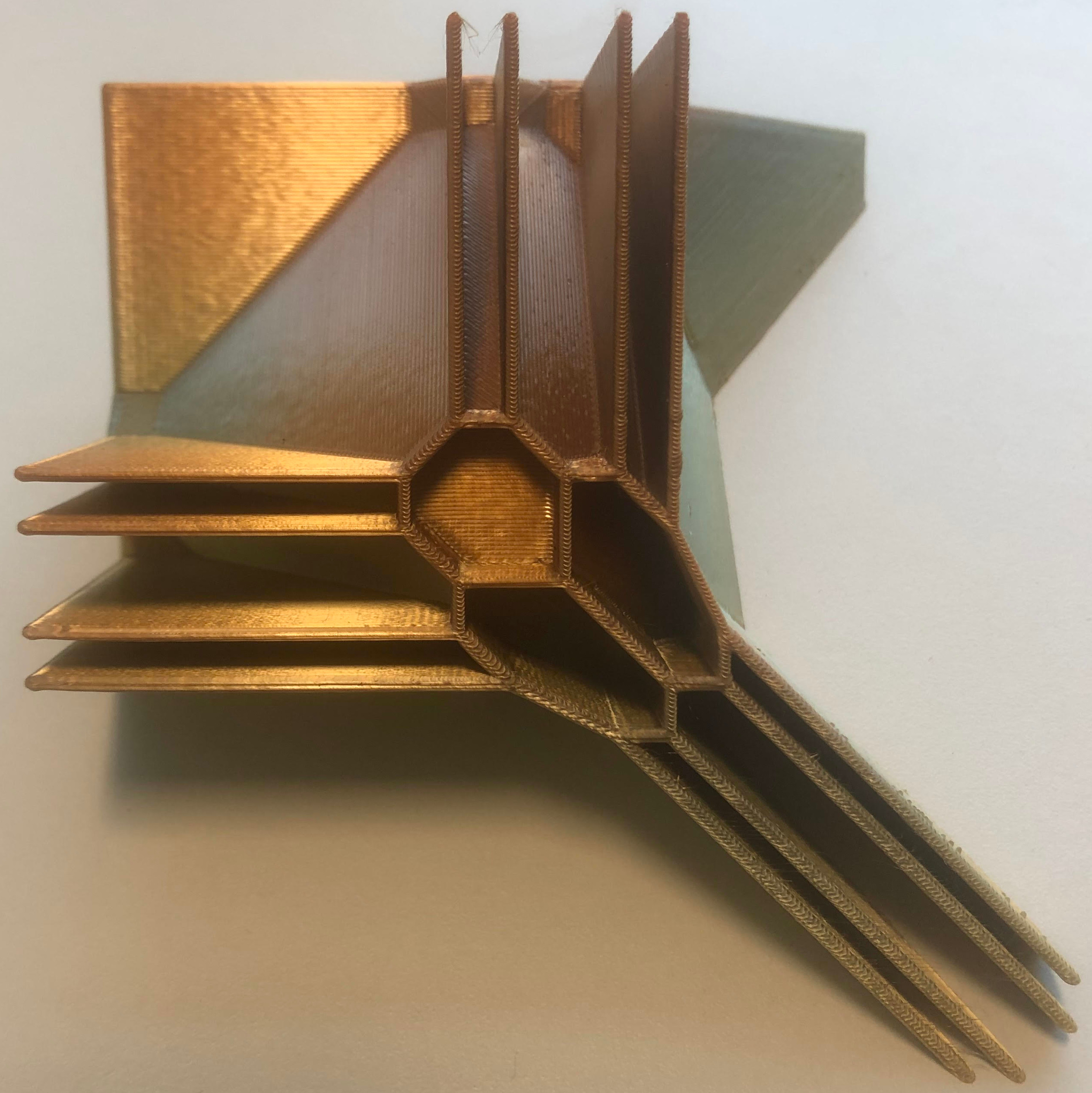}\hspace{5mm}%
 \includegraphics[height=37mm]{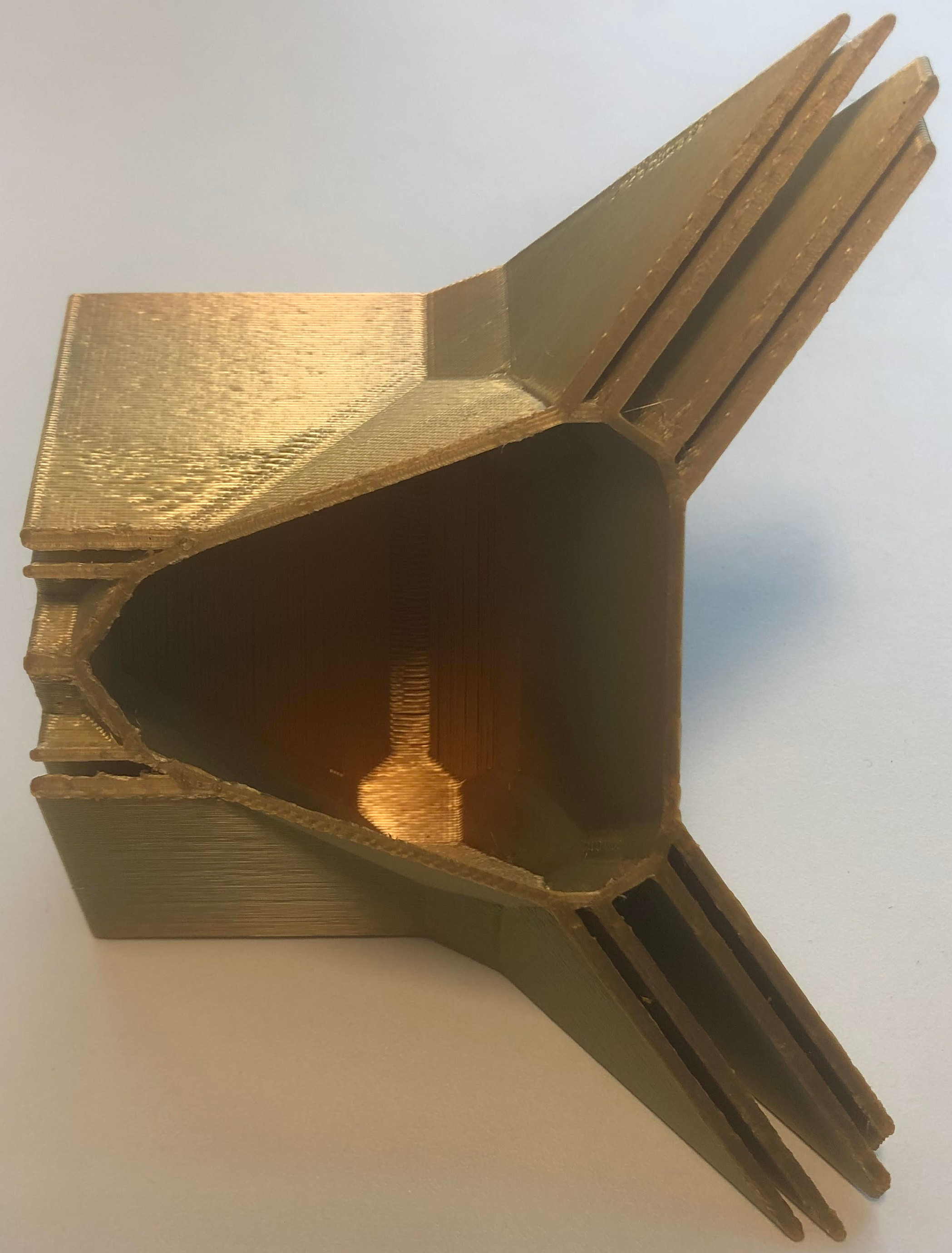}
 \caption{A tropical quartic surface from \cite{BPS21}, cut in half to reveal the inscribed K3 polytope.\label{fig:quarticSurface}}
\end{figurehere}

\begin{figurehere}
 \centering
 \includegraphics[width=0.8\columnwidth, height=62mm]{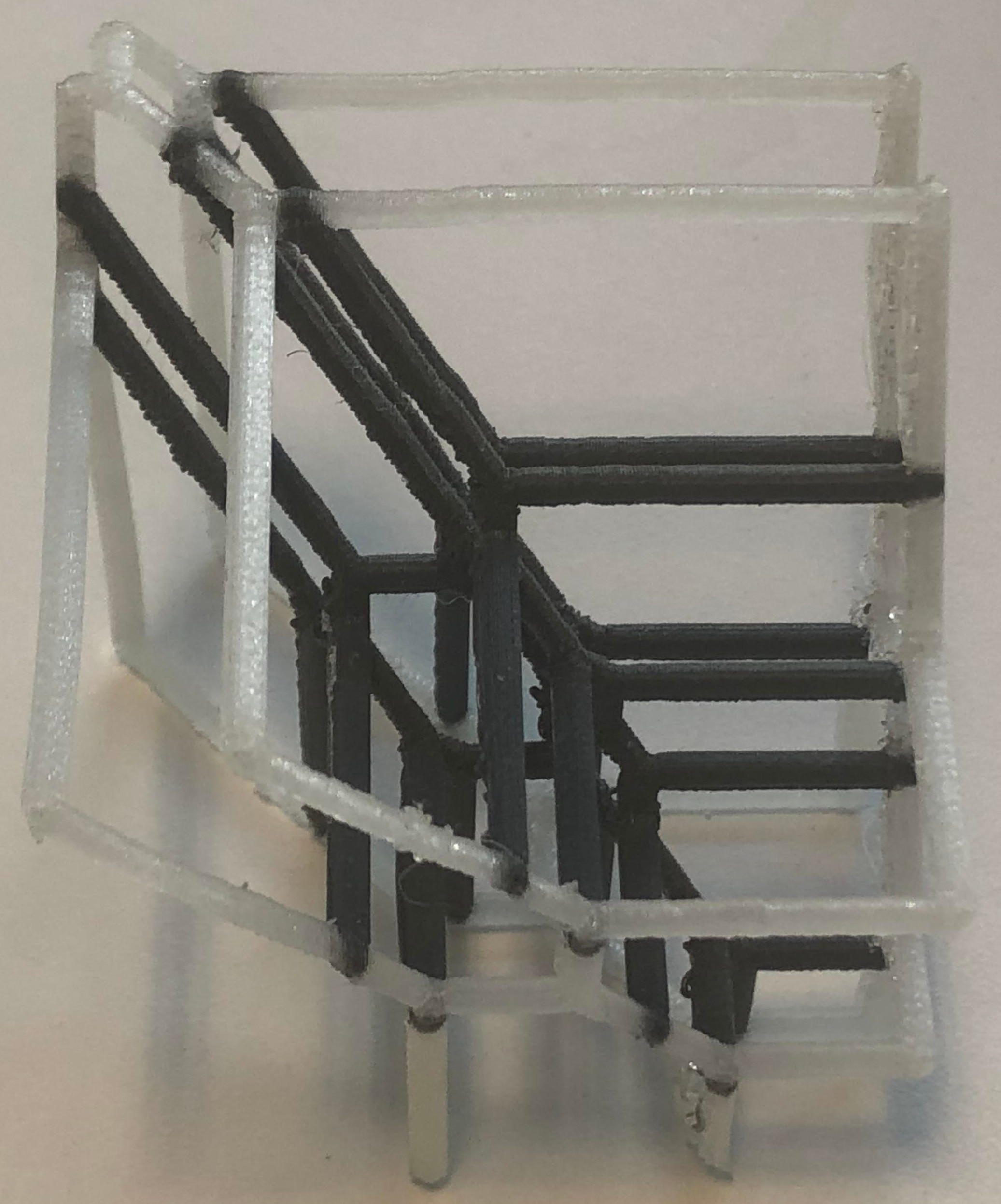}
 \caption{A tropical sextic curve from \cite{HL17}, the white frame traces the quadric surface containing it.\label{fig:sexticCurve}}
\end{figurehere}

%
%
%



\end{multicols}


\end{document}